\documentclass[12pt]{article}
\usepackage{latexsym,amsmath,amssymb}
\usepackage{color}
\definecolor{DarkBlue}{rgb}{0.1,0.1,0.5}
\definecolor{Red}{rgb}{0.9,0.0,0.1}

\newtheorem{theorem}{Theorem}[section]

\newtheorem{lemma}[theorem]{Lemma}

\def \BL{{\Bbb L}}
\def \BM{{\Bbb M}}
\def \BN{{\Bbb N}}

\def \BR{{\Bbb R}}

\def \qed{\hfill $\Box$ \vskip 5pt}

\begin{document}

\title{Moment convergence of $Z$-estimators 
and $Z$-process method for change point problems}
\author{
Ilia Negri \ \ and \ \ Yoichi Nishiyama\footnote{Corresponding author.}
\\ \\
{\small\em
Department of Information Technology and Mathematical Methods, University of Bergamo} \\
{\small\em
Viale Marconi, 5, 24044, Dalmine (BG), Italy}\\
{\small\tt ilia.negri@unibg.it}\\
and \\
{\small\em
The Institute of Statistical Mathematics} \\
{\small\em
10-3 Midori-cho, Tachikawa, Tokyo 190-8562, Japan} \\
{\small\tt nisiyama@ism.ac.jp}
}
\maketitle


\begin{abstract}
The problem to establish not only the asymptotic distribution results for statistical estimators 
but also the moment convergence of the estimators has been recognized as an important issue 
in advanced theories of statistics. 
One of the main goals of this paper is to present a metod to derive 
the moment convergence of $Z$-estimators as it has been done for $M$-estimators. 
Another goal of this paper is to develop a general, unified approach, based on 
some partial estimation functions which we call ``$Z$-process'', to the change point problems 
for ergodic models as well as some models where the Fisher information matrix is random 
and inhomogeneous in time. 
Applications to some diffusion process models and Cox's regression model are also discussed. 
\end{abstract}

%



\section{Introduction}
This paper is devoted to the study of two problems both based on 
the ``$Z$-methods'', in other words, some methods using 
the solutions to estimating equations. 
Let us first describe the outlines of the two themes, 
and next we shall list up some examples to which our results can be applied.  

\subsection{Theme I: Moment convergence of $Z$-estimators}
For an illustration, let us consider the simplest case of i.i.d.\ data. 
Let $ ({\cal X},{\cal A},\mu) $ be a measure space, 
and let us be given a parametric family of probability densities $ f(\cdot;\theta) $ 
with respect to $ \mu $, where $ \theta \in \Theta \subset \BR^d $. 
Let $ X_1,X_2,... $ be an i.i.d.\ sequence of $ {\cal X} $-valued random variables from this parametric model. 
There are two ways to define the ``maximum likelihood estimator (MLE)'' in statistics. 
One way is to define it as the maximum point of the random function 
\[
\theta \mapsto \BM_n(\theta)=\frac{1}{n}\sum_{k=1}^n \log f(X_k;\theta), 
\]
while the other is to do it as the solution to the estimating equation 
\[
\dot{\BM}_n(\theta)=0, 
\]
where $ \dot{\BM}_n(\theta) $ is the gradient vector of $ \BM_n(\theta) $. 
The former is a special case of ``$ M $-estimators'', and the latter is that 
of ``$ Z $-estimators''; see van der Vaart and Wellner (1996) for these terminologies. 
It may appear from the above introduction that $M$-estimation and $Z$-estimation can be regarded as ``almost equivalent''. 
However, it is not always true as we will discuss below. 

It is well known that the MLE $ \widehat{\theta}_n $ is asymptotically normal: 
it holds for any {\em bounded} continuous function $ \psi: \BR^d \to \BR $ that 
\[
\lim_{n \to \infty}
E[\psi(\sqrt{n}(\widehat{\theta}_n-\theta_0))] = E[\psi(I(\theta_0)^{-1/2}Z)], 
\]
where $ I(\theta_0) $ is the Fisher information matrix and $ Z $ is 
a standard Gaussian random vector. 
Furthermore, it is important for some advanced theories in statistics, including 
asymptotic expansions and model selections, 
to extend this kind of results for {\em bounded} continuous functions $ \psi $ 
to that for any continuous function $ \psi $ with {\em polynomial growth}, 
that is, any continuous function $ \psi $ for which there exist some constants $ C=C_\psi >0 $ and $ q=q_\psi >0 $ such that 
\begin{equation}\label{polynomial}
|\psi(x)|\leq C(1+||x||)^q, \quad \forall x \in \BR^d. 
\end{equation}
See the discussion in Yoshida (2011) for the importance of this problem. 

We observe that, when we have an asymptotic distribution result of an estimator, 
namely $ R_n(\widehat{\theta}_n-\theta_0) \to^d L(\theta_0) $ 
where $ R_n $ is a (possibly, random) diagonal matrix and the limit random vector $ L(\theta_0) $ is not necessarily Gaussian, 
it is sufficient for the generalization to the case where $ \psi $ is a continuous function 
satisfying (\ref{polynomial}) to check that $ ||R_n(\widehat{\theta}_n-\theta_0)|| $ is 
{\em asymptotically $ L_p $-bounded} for some $ p > q $, that is, 
\[
\limsup_{n \to \infty}E[||R_n(\widehat{\theta}_n-\theta_0)||^p] < \infty. 
\]

The study to provide some methods to obtain the moment convergence with polynomial order 
goes back to Ibragimov and Has'minskii (1981) 
who considered the MLEs and the Bayes estimators (as some special cases of $M$-estimators) 
in the general framework of the locally asymptotically normal models. 
It should be emphasized that one of the important merits of Ibragimov and Has'minskii's program is that 
the theory, based on the likelihood, automatically yields also the asymptotic efficiency. 
In their main theorems, it was assumed that an {\em exponential type large deviation inequality} 
holds for the rescaled log-likelihood ratio random field. 
Kutoyants (1984, 1994, 1998, 2004) sucessfully applied this theory to different stochastic process models 
using some characteristics of the models under consideration such as small diffusion models, ergodic diffusion models and Poisson process models. 
However, developing a general theory to establish the large deviation inequality had been an open problem for many years. 
This problem was solved, and the results have been published in Yoshida (2011). 
The paper starts from pointing out that a {\em polynomial type large deviation inequality} 
is sufficient for the core part of Igragimov and Has'minskii's (1981) program, 
and then the (polynomial type) large deviation inequality has been proved with a good generality. 
Uchida and Yoshida (2012) applied Yoshida's (2011) theory to establish the moment convergence of 
three kinds of adaptive $M$-estimators in ergodic diffusion process models including the one introduced by Kessler (1995). 
We mention that Nishiyama (2010) pointed out that the moment convergence problem 
for $M$-estimators can be solved also by using a maximal inequality instead of the large deviation inequalities, 
and that Kato (2011) took this type of approach to deal with some bootstrap $M$-estimators. 

In this paper, we consider the problem of the moment convergence of $Z$-estimators. 
Since we have to assume that the random field (something like the log-likelihood) is differentiable, 
the framework of $Z$-estimation is more restrictive than that of $M$-estimation. 
On the other hand, our proof is a combination of arguments based only on 
usual H\"older's and Minkowskii's inequalities, 
and no large deviation type inequality appears in our treatment for $Z$-estimators. 

Another difference between $M$- and $Z$-estimations is that in the latter theory 
the case where the rates of convergence are different over the components of $ \theta $ can be treated easily. 
This is due to the fact that in the theory of $Z$-estimation 
the gradient vector $ \dot{\BL}_n(\theta) $ of a contrast function $ \BL_n(\theta) $, where $ \BL_n(\theta) $ is typically the log-likelihood function, 
can by pre-multiplied by a {\em matrix} $ R_n^{-2} $ to get a kind of law of large numbers, namely, 
\[
\dot{\BM}_n(\theta)=R_n^{-2} \dot{\BL}_n(\theta). 
\]
Typically, $ R_n=\sqrt{n}I_d $ where  $ I_d $ is the identity matrix, 
although in the approach to $Z$-estimation presented here the diagonal components of $ R_n $ may be different. 
Compare this with the framework of $M$-estimation where the (scalar valued) contrast function $ \BL_n(\theta) $ 
with no assumption of differentiability has to be multiplied by a scalar. 
Yoshida (2011) dealt with this point developing an iterative method. 

Thinking of these differences, we may conclude that $M$-estimation and $Z$-estimation are not ``almost equivalent'' 
at least for the moment convergence problem. 

\subsection{Theme II: $Z$-process method for change point problems} 
Let us give an illustration by the example of independent data again. 
We introduce the partial sum process  
\[
\BM_n(u,\theta)=\frac{1}{n}\sum_{k=1}^{[un]} \log f(X_k;\theta), \quad \forall u \in [0,1], 
\]
and consider the gradient vectors $ \dot{\BM}_n(u,\theta) $ of $ \BM_n(u,\theta) $ with respect to $ \theta$. 
Let $ \widehat{\theta}_n $ be the MLE for the full data $ X_1,...,X_n $ 
as a special case of $Z$-estimators, that is, $ \widehat{\theta}_n $ is the solution to the estimating equation 
\[
\dot{\BM}_n(1,\theta)=0. 
\]

The fact that the random process 
\[
u \leadsto \sqrt{n}\dot{\BM}_n(u,\theta_0) 
\quad \mbox{converges weakly to } \quad 
u \leadsto I(\theta_0)^{1/2}B(u)
\]
in the Skorohod space $ D[0,1] $, where $ u \leadsto B(u) $ 
is a vector of independent standard Brownian motions, is immediate 
from Donsker's theorem. 
However, 
it does not seem so well known that the random process 
\begin{equation}\label{not monkey}
u \leadsto \sqrt{n}\dot{\BM}_n(u,\widehat{\theta}_n) 
\quad \mbox{converges weakly to} \quad 
u \leadsto I(\theta_0)^{1/2}B^\circ(u)
\end{equation}
in $ D[0,1] $, where $ u \leadsto B^\circ(u) $ 
is a vector of independent standard Brownian bridges. 
Horv\'{a}th and Parzen (1994) is apparently the first to introduce the statistic 
\[
{\cal T}_n= n \sup_{u \in [0,1]}\dot{\BM}_n(u,\widehat{\theta}_n)^\top \widehat{I}_n^{-1}\dot{\BM}_n(u,\widehat{\theta}_n) 
\]
for change point problems, where $ \widehat{I}_n $ is a consistent estimator for the Fisher Information matrix $ I(\theta_0) $. 
It is immediate from (\ref{not monkey}) and the continuous mapping theorem that 
\[
{\cal T}_n \to^d \sup_{u \in [0,1]}||B^\circ(u)||^2. 
\]
Let us call this approach pioneered by Horv\'{a}th and Parzen (1994) ``$Z$-process method''. 

Although Horv\'{a}th and Parzen (1994) didn't discuss the asymptotic behavior of the test 
under the alternative, Negri and Nishiyama (2011) who took the $Z$-process method 
for an ergodic diffusion process model based on the continuous observation 
proved also the consistency of the test under an alternative which has sufficient generality. 
Negri and Nishiyama's (2011) argument for alternatives can be applied also to the case of independent data. 
In Section \ref{II} of this paper, we will present a generalized version of $Z$-process method which works also for some cases where 
the limit of the test statistic under the null hypothesis is a functional of 
a ``mixture'' of standard Brownian motions, 
where the mixture process, which is something like a ``partial process of Fisher information'', is random and time dependent. 
We will also develop the argument of Negri and Nishiyama (2011) for the consistency under alternatives in a more general way. 
Some new examples will be given. 

\subsection{Notations and examples}\label{introduction of examples}
In the rest of this section, we shall list up some examples to which our results can be applied. 
In what follows, the parameter space $ \Theta $ is a bounded, open, convex subset of $ \BR^d $, 
where $ d $ is a fixed, positive integer. 
The word ``vector'' always means ``$ d $-dimensional real column vector'', 
and the word ``matrix'' does ``$ d \times d $ real matrix''. 
The Euclidean norm is denoted by $ ||v||:=\sqrt{\sum_{i=1}^d |v^{(i)}|^2} $ 
for a vector $v$ where $ v^{(i)} $ denotes the $ i $-th component of $v $, 
and by $ ||A||:=\sqrt{\sum_{i,j=1}^d|A^{(i,j)}|^2} $ for a matrix $A$ 
where $ A^{(i,j)} $ denotes the $ (i,j) $-component of $ A $. 
Note that $ ||Av|| \leq ||A||\cdot ||v|| $ and $ ||AB|| \leq ||A||\cdot||B|| $ for 
vector $ v $ and matrices $ A, B $. 
The notations $ v^\top $ and $ A^\top $ denote the transpose. 
We use also the notation $ A \circ B $ defined by 
$ (A \circ B)^{(i,j)}:=A^{(i,j)}B^{(i,j)} $ for two matrices $ A, B $ (the Hadamard product). 
We denote by $ I_d $ the identity matrix. 
The notations $ \to^p $ and $ \to^d $ mean the convergence in probability and 
the convergence in distribution, as $ n \to \infty $, respectively. 

\vskip 10pt
{\em Example A: Ergodic diffusion process.} 
Let $ I = (l,r) $, where $ -\infty \leq l < r \leq \infty $, be given. 
Let us consider an $ I $-valued diffusion process $ t \leadsto X_t $ which 
is the unique strong solution to the stochastic differential equation (SDE) 
\[
X_t=X_0+ \int_0^t S(X_s;\alpha)ds + \int_0^t \sigma(X_s;\beta)dW_s, 
\]
where $ s \leadsto W_s $ is a standard Wiener process. 
The parameters come from $ \alpha \in \Theta_A \subset \BR^{d_A} $ 
and $ \beta \in \Theta_B \subset \BR^{d_B} $, 
and we denote $ \theta= (\alpha^{\top}, \beta^{\top})^{\top} $. 
We are supposed to be able to observe the process $ X $ 
at discrete time grids $ 0=t_0^n< t_{1}^n < \cdots < t_n^n $, 
and we shall consider the asymptotic scheme $ n\Delta_n^2 \to 0 $ 
and $ t_n^n \to \infty $ as $ n \to \infty $, where 
\[
\Delta_n=\max_{1 \leq k \leq n}|t_{k}^n-t_{k-1}^n|, 
\]
and 
\begin{equation}\label{equidistant}
\sum_{k=1}^n \left|\frac{|t_k^n -t_{k-1}^n|}{t_n^n} - \frac{1}{n} \right|
\to 0, \quad \mbox{as } n \to \infty. 
\end{equation}

For Themes I and II, we introduce 
\[
\dot{\BM}_n(\theta)=R_n^{-2} \dot{\BL}_n(1,\theta) \quad \mbox{and} \quad 
\dot{\BM}_n(u,\theta)=R_n^{-2} \dot{\BL}_n(u,\theta),
\]
respectively, where 
\begin{eqnarray*}
\lefteqn{
\BL_n(u,\theta)}
\\
&=&-\sum_{k:  t_{k-1}^n \leq u t_n^n}\left\{ \log \sigma(X_{t_{k-1}^n};\beta) 
+\frac{|X_{t_{k}^n}-X_{t_{k-1}^n}-S(X_{t_{k-1}^n};\alpha)|t_k^n - t_{k-1}^n||^2}
{2\sigma(X_{t_{k-1}^n};\beta)^2|t_{k}^n-t_{k-1}^n|} \right\} 
\end{eqnarray*}
and $ R_n $ is the diagonal matrix such that $ R_n^{(i,i)} $ 
is $ \sqrt{t_n^n} $ for $ i=1,...,d_A $ and $ \sqrt{n} $ for $ i=d_A+1,...,d $ 
with $ d=d_A+d_B $. 

The problem to establish 
the moment convergence for $ M $-estimators in this model, 
where $ X $ is a multi-dimensional diffusion process, was considered by Yoshida (2011), 
and Uchida and Yoshida (2012) relaxed the assumption $ n\Delta_n^2 \to 0 $ up to $ n \Delta_n^{a}\to 0 $, where $ a \geq 2 $ 
is a constant depending on the smoothness of the model as was done by Kessler (1995, 1997). 
To treat the parameters with different rates of convergence, in both papers an iterative method is used, and 
the method leads to some ``adaptive estimators'' that have advantages in applications as their simulation results show. 
In order to explain our core idea clearly, we only consider the one-dimensional diffusion process $ X $ under the sampling scheme $ n\Delta_n^2 \to 0 $. 
Some extension with no interative argument  to the case considered in Uchida and Yoshida (2012) could be possible. 

Regarding Theme 2, Song and Lee (2009) proposed a statistic for testing the existence of a change point 
of the parameter $ \beta $, but the problem to test it for both parameters 
was left as an open problem in their paper (see their Section 5). 
We will give an answer to this problem in Section \ref{Song Lee}. 

\vskip 10pt
{\em Example B: Volatility of diffusion process.} 
Let $ I = (l,r) $, where $ -\infty \leq l < r \leq \infty $, be given.  
Let us consider an $ I $-valued diffusion process $ t \leadsto X_t $ which 
is the unique strong solution to the SDE 
\[
X_t=X_0+ \int_0^t S(X_s)ds + \int_0^t \sigma(X_s;\theta)dW_s, 
\]
where $ s \leadsto W_s $ is a standard Wiener process. 
Here, the drift coefficient $S(\cdot)$ is treated as an unknown nuisance function. 
We are supposed to be able to observe the process $ X $ 
at discrete time grids $ 0=t_0^n< t_{1}^n < \cdots < t_n^n=T< \infty $, 
and we shall consider the asymptotic scheme (\ref{equidistant}). 

We introduce 
\[
\dot{\BM}_n(\theta)=\frac{1}{n} \dot{\BL}_n(1,\theta) \quad 
\mbox{and} \quad \dot{\BM}_n(u,\theta) =\frac{1}{n}\dot{\BL}_n(u,\theta) 
\]
for Themes I and II, respectively, where 
\[
\BL_n(u,\theta)=-\sum_{k:  t_{k-1}^n \leq u t_n^n}\left\{ \log \sigma(X_{t_{k-1}^n};\theta) 
+\frac{|X_{t_{k}^n}-X_{t_{k-1}^n}|^2}
{2\sigma(X_{t_{k-1}^n};\theta)^2|t_{k}^n-t_{k-1}^n|} \right\}. 
\]
The rate matrix is given by $ R_{n}=\sqrt{n}I_d $. 

Iacus and Yoshida (2012) proposed an estimator for the change point in a similar model. 
Our result for Theme II here, dealing with testing the existence of a change point, can be applied before statisticians proceed to their theory of estimation. 

As we already mentioned, Song and Lee (2009) proposed a statistic for testing the existence of 
a change point in the volatility of an ergodic diffusion process 
under the asymptotic scheme $ n\Delta_n^a \to 0 $ and $ n\Delta_n^b \to \infty $ for some $ a > b > 4 $, 
by an approach which is different from ours. 

\vskip 10pt
{\em Example C: Cox's regression model.} 
Let a sequence of counting processes $ t \leadsto N_t^k $, $ k=1,2,... $, 
which do not have simultaneous jumps, be observed during the time interval $ [0,T] $. 
Suppose that $ t \leadsto N_t^k $ has the intensity 
\[
\lambda_t^k(\theta)=\alpha(t) e^{\theta^\top Z_t^k} Y_t^k, 
\]
where the baseline hazard function $ \alpha $ which is common for all $ k $'s is non-negative and 
satisfies $ \int_0^T\alpha(t)dt< \infty $, 
the random process $ t \leadsto Z_t^k $ is a vector valued covariate for the individual $ k $, and 
the random process $ t \leadsto Y_t^k $ is given by 
\[
Y_t^k=\left\{
\begin{array}{lll}
1, & & \mbox{if the individual } k \mbox{ is observed at time $t$}, 
\\
0, & & \mbox{otherwise}. 
\end{array}
\right.
\]
This model was introduced by Cox (1972), and its asymptotic theory was developed by 
Andersen and Gill (1982). 

For Themes I and II, we introduce 
\[
\dot{\BM}_n(\theta)=\frac{1}{n} \dot{\BL}_n(1,\theta) 
\quad \mbox{and} \quad \dot{\BM}_n(u,\theta)=\frac{1}{n}\dot{\BL}_n(u,\theta), 
\]
respectively, where 
\[
\BL_n(u,\theta)=\sum_{k=1}^{n}\int_0^{uT} (\theta^\top Z_t^k 
- \log S_t^{n,0}(\theta) )dN_t^k 
\]
with 
\[
S_t^{n,0}(\theta)=\sum_{k=1}^n e^{\theta^\top Z_t^k} Y_t^k. 
\]
The rate matrix is $ R_n=\sqrt{n} I_d $. 

\vskip 10pt
{\em Example D: Counting process models.} 
Let $ t \leadsto N_t $ be a counting process with the intensity $ t \leadsto \lambda_t(\theta) $. 
We suppose that we can observe the processes on the compact time interval $ [0,T_n] $, 
and consider the asymptotic scheme $ T_n \to \infty $. 
Our results may be typically applied to 
\[
\dot{\BM}_n(\theta)=\frac{1}{T_n}\dot{\BL}_n(1,\theta) \quad \mbox{and} \quad 
\dot{\BM}_n(u,\theta)=\frac{1}{T_n}\dot{\BL}_n(u,\theta), 
\]
where 
\[
\BL_n(u,\theta)=
\int_0^{uT_n} \log \lambda_t(\theta)dN_t -\int_0^{u T_n} \lambda_t(\theta)dt. 
\]
The rate matrix is typically $ R_n=\sqrt{T_n}I_d $. 
An example of this model is the stress release process introduced by Isham and Westcott (1979): 
\[
\lambda_t(\alpha,\beta,\theta)=\phi(\alpha t-\beta N_{t-};\theta). 
\]
It is known that the process $ t \leadsto X_t:=\alpha t-\beta N_{t-} $ is ergodic 
under some mild conditions. 
A test statistic for the change point problem in this model, which is different from ours, 
has been proposed by Fujii and Nishiyama (2011).   

\vskip 10pt
{\em Example E: Non-linear time series models.} 
Let us consider the time series models of the form 
\[
X_k=a(X_{k-1},X_{k-2},...;\theta) + b(X_{k-1},X_{k-2},...;\theta)\varepsilon_k, 
\quad k=1,2,.... 
\]
Here, $ \{ \varepsilon_k \} $ is an i.i.d.\ sequence with $ E[\varepsilon_1]=0 $ or, 
more generally, a martingale difference sequence with respect to the 
filtration $ ({\cal F}_k)_{k \geq 0} $ where $ {\cal F}_k=\sigma(X_k,X_{k-1},...) $. 
A possible way to define the estimating functions is 
\[
\dot{\BM}_n(\theta)=\frac{1}{n}\dot{\BL}_n(1,\theta) 
\quad \mbox{and}\quad 
\dot{\BM}_n(u,\theta)=\frac{1}{n}\dot{\BL}_n(u,\theta), 
\]
where  
\[
\BL_n(u,\theta)=-\sum_{k=1}^{[un]}\left\{\log b(X_{k-1},X_{k-2},...;\theta) 
+\frac{|X_{k}-a(X_{k-1},X_{k-2},...;\theta)|^2}
{2b(X_{k-1},X_{k-2},...;\theta)^2} \right\}. 
\]
The rate matrix is typically given by $ R_{n}=\sqrt{n}I_d $. 

\vskip 10pt
Some detailed discussions for Examples A, B, and C will be given in Sections \ref{A}, \ref{B} and \ref{C}, respectively, 
while Examples D and E are not discussed in detail here. 

\section{Moment convergence of $Z$-estimators}\label{I}

Let $ \Theta $ be a bounded, open, convex subset of $ \BR^d $. 
Let us be given a real valued random function $ \BL_n(\theta) $ of $ \theta \in \Theta $ 
which is twice continuously differentiable 
with the gradient vector $  \dot{\BL}_n(\theta) $ and the Hessian matrix $ \ddot{\BL}_n(\theta) $, defined on a probability space 
$ (\Omega,{\cal F},P) $ that is common for all $ n \in \BN $. 
(However, it will be clear from our proofs that if the limit matrices $ V(\theta_0) $ and $ \ddot{M}(\theta) $ 
appearing below are non-random then the underlying probability spaces need not to be common for 
all $ n \in \BN $.) 
Let $ R_n $ be a (possibly, random) diagonal matrix whose diagonal components are positive, 
and define $ Q_n $ by $ Q_n^{(i,j)}=(R_{n}^{(i,i)} R_{n}^{(j,j)})^{-1} $. Using these matrices, put 
\begin{equation}\label{definition of gradient and Hessian}
\dot{\BM}_n(\theta) :=R_n^{-2}\dot{\BL}_n(\theta) \quad \mbox{and} \quad
\ddot{\BM}_n(\theta) :=Q_n \circ \ddot{\BL}_n(\theta). 
\end{equation}
(In the typical cases, $ R_n=\sqrt{n}I_d $ and $ Q_n=n^{-1} {\bf 1} $, where $ {\bf 1} $ denotes the matrix whose all 
components are $1$.) 

First, we state a theorem to give an asymptotic representation for $Z$-estimators. 
Although this result is not really novel, we will give a proof for references. 

\begin{theorem}\label{main theorem 1-1}
Consider the above setting. 
Suppose that there exists a sequence of matrices $ V_n(\theta_0) $ 
which are regular almost surely 
such that for any sequence of $ \Theta $-valued random vectors 
$ \widetilde{\theta}_n $ converging in probability to $ \theta_0 $, 
\[
\ddot{\BM}_n(\widetilde{\theta}_n)+V_n(\theta_0) \to^p 0. 
\]
Suppose also that 
\[
(R_n\dot{\BM}_n(\theta_0),V_n(\theta_0)) \to^d (L(\theta_0), V(\theta_0)), 
\]
where $ L(\theta_0) $ is a random vector and 
$ V(\theta_0) $ is a random matrix which is regular almost surely 
(we do not assume that $ V(\theta_0) $ and $ L(\theta_0) $ are independent). 

Then, for any sequence of $ \Theta $-valued random vectors 
$ \widehat{\theta}_n $ which converges in probability to $ \theta_0 $ 
and satisfies $ ||R_n\dot{\BM}_n(\widehat{\theta}_n)||=o_{P}(1) $, 
it holds that 
\begin{eqnarray*}
R_n(\widehat{\theta}_n-\theta_0) &=& V_n(\theta_0)^{-1}R_n\dot{\BM}_n(\theta_0)
+o_{P}(1)\\
&\to^d& 
V(\theta_0)^{-1}L(\theta_0). 
\end{eqnarray*}
\end{theorem}

{\bf Remark.} 
In this theorem, the consistency of the sequence of $Z$-estimators $ \widehat{\theta}_n $ 
has been assumed. A method to show this property will be 
given in Lemma \ref{lemma consistency} (i) below. 

Now, we give a theorem to establish the moment convergence of $Z$-estimators, 
which is the main result of this section. 

\begin{theorem}\label{main theorem 1-2}
Consider the setting described in the first paragraph of this section. 
Let some constants $ p \geq 1 $ and $ a,b> 1 $ such that $ \frac{1}{a}+\frac{1}{b}=1 $ be given; 
see a remark at the end of the theorem for the case where we may set $ a=1 $. 

Suppose that 
\begin{equation}\label{main theorem 1-2 a condition}
||R_n\dot{\BM}_n(\theta_0)|| \ \ \mbox{is asymptotically } L_{pa} \mbox{-bounded}. 
\end{equation}
Suppose also that there exist a constant $ \gamma \in (0,1] $ and some random matrices $ \ddot{M}(\theta) $ 
indexed by $ \theta \in \Theta $ such that 
\begin{equation}\label{main theorem 1-2 key condition}
\lim_{n \to \infty}E\left[\sup_{\theta \in \Theta}||R_n^{\gamma} (\ddot{\BM}_n(\theta)-\ddot{M}(\theta))||^{pa/\gamma}\right] = 0. 
\end{equation}
Suppose further that either of the following [M1] or [M2] is satisfied: 

[M1] There exists a random matrix $ J(\theta_0) $ which is positive definite almost surely such that 
$ \ddot{M}(\theta) \leq -J(\theta_0) $ for all $ \theta \in \Theta $, almost surely, 
and that $ E[||J(\theta_0)^{-1}||^{pb/\gamma}] < \infty  $; 

[M2] $ E[\sup_{\theta \in \Theta}||\ddot{M}(\theta)^{-1}||^{pb/\gamma}] < \infty $, 
where the random matrices $ \ddot{M}(\theta) $'s are assumed to be regular almost surely. 
 
Then, for any sequence of $ \Theta $-valued random vectors $ \widehat{\theta}_n $ 
such that $ ||R_n\dot{\BM}_n(\widehat{\theta}_n)|| $ is asymptotically $ L_{pa} $-bounded, 
it holds that $ ||R_n(\widehat{\theta}_n-\theta_0)|| $ is asymptotically 
$ L_p $-bounded. 
Therefore, in this situation, whenever we also have that $ R_n(\widehat{\theta}_n-\theta_0) \to^d G(\theta_0) $ 
where $ G(\theta_0) $ is a random vector, 
it holds  for any continuous function $ \psi $ satisfying (\ref{polynomial}) for $ q \in (0,p) $ that 
\[
\lim_{n \to \infty}E[\psi(R_n(\widehat{\theta}_n-\theta_0))]
= E[\psi(G(\theta_0))], 
\]
where the limit is also finite. 

When the last condition in [M1] is satisfied with $ ||J(\theta_0)||^{-1} $ which is bounded or the first condition in 
[M2] is satisfied with 
$ \sup_{\theta \in \Theta}||\ddot{M}(\theta)^{-1}|| $ which is bounded, the constant $ a $ appearing in the above 
claim may be replaced by $ 1 $. 
\end{theorem}

{\bf Remark.} The condition [M1] corresponds to the case $ \rho=2 $ of 
the conditions [A3] and [A5] in Yoshida (2011), which are 
\[
M(\theta)-M(\theta_0) \leq -\chi (\theta_0)||\theta-\theta_0||^\rho, \quad \forall \theta \in \Theta, 
\]
and high order moment conditions on the positive random variable $ \chi(\theta_0)^{-1} $, 
where ``$M(\theta)$'' should be read as ``$Y(\theta)$'' in Yoshida's (2011) notation. 

\vskip 10pt
{\it Proof of Theorem \ref{main theorem 1-1}.} 
Recalling (\ref{definition of gradient and Hessian}), 
it follows from the Taylor expansion that 
\begin{equation}\label{Taylor}
(R_n\dot{\BM}_n(\widehat{\theta}_n))^{(i)}=(R_n\dot{\BM}_n(\theta_0))^{(i)}
+ (\ddot{\BM}_n(\widetilde{\theta}_n)R_n(\widehat{\theta}_n-\theta_0))^{(i)}, \quad i=1,...,d, 
\end{equation}
where $ \widetilde{\theta}_n $ is a random vector on the segment connecting $ \theta_0 $ and $ \widehat{\theta}_n $. 
So we have 
\begin{equation}\label{main theorem 1 zenhan key}
R_n(\widehat{\theta}_n-\theta_0)
=A_n + B_nR_n(\widehat{\theta}_n-\theta_0), 
\end{equation}
where 
\begin{eqnarray*}
A_n&=&V_n(\theta_0)^{-1}
R_n(\dot{\BM}_n(\theta_0)-\dot{\BM}_n(\widehat{\theta}_n)), 
\\
B_n&=&V_n(\theta_0)^{-1}(\ddot{\BM}_n(\widetilde{\theta}_n)+V_n(\theta_0)). 
\end{eqnarray*}
It follows from the extended continuous mapping theorem 
(e.g., Theorem 1.11.1 of van der Vaart and Wellner (1996)) 
that $ V_n(\theta_0)^{-1} \to^p V(\theta_0)^{-1} $, thus 
we have $ ||A_n||=O_{P}(1) $ and $ ||B_n||=o_{P}(1) $. 
It therefore holds that 
\[
||R_n(\widehat{\theta}_n-\theta_0)||\leq O_{P}(1) + o_{P}(1) \cdot ||R_n(\widehat{\theta}_n-\theta_0)||, 
\]
which implies that $ ||R_n(\widehat{\theta}_n-\theta_0)||=O_{P}(1) $. 
Hence, going back to (\ref{main theorem 1 zenhan key}) we obtain 
\[
R_n(\widehat{\theta}_n-\theta_0)=A_n+o_{P}(1)=V_n(\theta_0)^{-1}R_n\dot{\BM}_n(\theta_0)+ o_{P}(1). 
\]
The last claim is also a consequence of the extended 
continuous mapping theorem. 
The proof is finished. 
\qed

\vskip 10pt
{\it Proof of Theorem \ref{main theorem 1-2}.} 
We will give a proof for the case where [M1] is assumed. 
The proof for the case where [M2] is assumed is similar (and simpler), so it is omitted. 

Due to (\ref{Taylor}) again, we have  
\[
R_n(\widehat{\theta}_n-\theta_0)
=C_n + (D_n^{(1)}+D_n^{(2)})R_n(\widehat{\theta}_n-\theta_0), 
\]
where 
\begin{eqnarray*}
C_n&=&J(\theta_0)^{-1}
R_n(\dot{\BM}_n(\theta_0)-\dot{\BM}_n(\widehat{\theta}_n)), 
\\
D_n^{(1)} &=&J(\theta_0)^{-1}(\ddot{\BM}_n(\widetilde{\theta}_n)-\ddot{M}(\widetilde{\theta}_n)),
\\
D_n^{(2)} &=&J(\theta_0)^{-1}(\ddot{M}(\widetilde{\theta}_n)+J(\theta_0)),
\end{eqnarray*}
where $ \widetilde{\theta}_n $ is a random vector on the segment connecting $ \theta_0 $ and $ \widehat{\theta}_n $. 

From now on, we consider the case $ \gamma \in (0,1) $; the proof for the case $ \gamma=1 $ is easier, and 
it is omitted. 
Since $ -D_n^{(2)} $ is non-negative definite almost surely, 
it follows from Minkowski's and H\"older's inequalities that 
\begin{eqnarray*}
\lefteqn{
(E[||R_n(\widehat{\theta}_n-\theta_0)||^p])^{1/p}
}\\
&\leq&
(E[||(I_d-D_n^{(2)})R_n(\widehat{\theta}_n-\theta_0)||^p])^{1/p}
\\
&\leq&
(E[||C_n||^p])^{1/p}
+
(E[||R_n^{\gamma}D_n^{(1)}||^{p/\gamma}])^{\gamma/p}(E[||R_n^{1-\gamma}(\widehat{\theta}_n-\theta_0)||
^{p/(1-\gamma)}])^{(1-\gamma)/p}
\\
&\leq&
O(1) + o(1)\cdot (E[||R_n^{1-\gamma}(\widehat{\theta}_n-\theta_0)||^{p/(1-\gamma)}])^{(1-\gamma)/p}, 
\end{eqnarray*}
where we have used H\"older's inequality again to get 
\[
E[||C_n||^p] \leq (E[||J(\theta_0)^{-1}||^{pb}])^{1/b}(E[||R_n(\dot{\BM}_n(\theta_0)-\dot{\BM}_n(\widehat{\theta}_n))||^{pa})^{1/a}
\]
and 
\[
E[||R_n^{\gamma}D_n^{(1)}||^{p/\gamma}] \leq (E[||J(\theta_0)^{-1}||^{pb/\gamma}])^{1/b}
(E[||R^{\gamma}(\ddot{\BM}_n(\widetilde{\theta}_n)-\ddot{M}(\widetilde{\theta}_n))||^{pa/\gamma})^{1/a}; 
\]
if $ ||J(\theta_0)||^{-1} $ is bounded, we can get this kind of bounds with $ a=1 $. 

Notice that 
\begin{eqnarray*}
\lefteqn{
||R_n^{1-\gamma}(\widehat{\theta}_n-\theta_0)||^{1/(1-\gamma)}
}\\
&\leq &
\sqrt{d^{(1/(1-\gamma))-1}\sum_{i=1}^d |R_n^{(i,i)}|^2|\widehat{\theta}_n^{(i)}-\theta_0^{(i)}|^{2/(1-\gamma)}}
\\
&\leq& ||R_n(\widehat{\theta}_n-\theta_0)||\cdot d^{1/(2-2\gamma)} \cdot |{\cal D}(\Theta)|^{\gamma/(1-\gamma)}, 
\end{eqnarray*}
where $ {\cal D}(\Theta) $ denotes the diameter of $ \Theta $. 
So we obtain 
\begin{eqnarray*}
\lefteqn{
(E[||R_n(\widehat{\theta}_n-\theta_0)||^p])^{1/p}
}\\
&\leq&
O(1)+
o(1) \cdot 
(E[||R_n(\widehat{\theta}_n-\theta_0)||^{p}])^{(1-\gamma)/p}
\\
&\leq&
O(1)+
o(1) \cdot 
(E[||R_n(\widehat{\theta}_n-\theta_0)||^{p}] \vee 1)^{1/p}, 
\end{eqnarray*}
which yields that 
\[
E[||R_n(\widehat{\theta}_n-\theta_0)||^p]
\leq O(1)+ o(1) \cdot E[||R_n (\widehat{\theta}_n-\theta_0)||^{p}]. 
\]
Therefore, $ ||R_n(\widehat{\theta}_n-\theta_0)|| $ is asymptotically $ L_p $-bounded. 
\qed

\section{$Z$-process method for change point problems}\label{II}

Let $ D[0,1] $ be the space of functions defined on $ [0,1] $ 
taking values in a finite-dimensional Euclidean space, 
which are right continuous and have left hand limits; 
we equip this space with the Skorohod metric. 
Throughout this section, all random processes, 
denoted as $ u \leadsto X(u) $, are assumed to take values in $ D[0,1] $. 

Let $ \Theta $ be a bounded, open, convex subset of $ \BR^d $. 
For every $ n \in \BN $, let $ u \leadsto \BL_n(u,\theta) $ 
be a real valued random process indexed by $ \theta \in \Theta $, 
defined on a probability space 
$ (\Omega,{\cal F},P) $ that is common for all $ n \in \BN $. 
(However, it will be clear from our proofs that 
the underlying probability spaces do not have to be common 
for $ n \in \BN $ 
if the objects $ \dot{M}_{\theta_0}(u,\theta) $, $ \dot{{\cal M}}(u,\theta) $ and $ V(u,\theta_0) $ appearing 
in the limit below are non-random.) 
We suppose that for every $ u \in [0,1] $ the random function 
$ \theta \mapsto \BL_n(u,\theta) $ is two times continuously differentiable 
with the gradient vector $ \dot{\BL}_n(u,\theta) $ and the Hessian matrix $ \ddot{\BL}_n(u,\theta) $. 
Let a (possibly, random) diagonal matrix $ R_n $ whose diagonal components are positive be given, 
and define $ Q_n $ by $ Q_n^{(i,j)}=(R_n^{(i,i)}R_n^{(j,j)})^{-1} $. Using these matrices, put 
\[
\dot{\BM}_n(u,\theta):=R_n^{-2} \dot{\BL}_n(u,\theta) \quad \mbox{and} \quad
\quad \ddot{\BM}_n(u,\theta) := Q_n \circ \ddot{\BL}_n(u,\theta).  
\]

We consider the following testing problem: 

\vskip 5pt
$ H_0 $: the true value $ \theta_0 \in \Theta $ 
does not change during $ u \in [0,1] $; 

$ H_1 $: ``not $ H_0 $''. 

\vskip 5pt
The meaning of ``not $ H_0 $'' will be precisely specified in the condition [A] below. 
Let us describe some properties which the ``limits'' $ \dot{M}_{\theta_0}(u,\theta) $ and $ \dot{{\cal M}}(u,\theta) $ of the random vectors 
$ \dot{\BM}_n(u,\theta) $ under $ H_0 $ and under $ H_1 $, respectively, have to satisfy. 

\vskip 5pt
[N] Under $ H_0 $, it holds that 
\begin{equation}\label{null 1}
\sup_{\theta \in \Theta}
||\dot{\BM}_n(1,\theta)- \dot{M}_{\theta_0}(1,\theta)||\to^p 0, 
\end{equation}
where the limits $ \dot{M}_{\theta_0}(1,\theta) $'s satisfy that 
\begin{equation}\label{null 2}
\inf_{\theta: ||\theta-\theta_0||> \varepsilon}
||\dot{M}_{\theta_0}(1,\theta)||
>0=||\dot{M}_{\theta_0}(1,\theta_0)||, 
\quad \mbox{almost surely}, \quad \forall \varepsilon > 0. 
\end{equation}

[A] Under $ H_1 $, it holds that  
\begin{equation}\label{alternative 1}
\sup_{u \in [0,1]}\sup_{\theta \in \Theta}||\dot{\BM}_n(u,\theta)-\dot{{\cal M}}(u,\theta)|| \to^p 0, 
\end{equation}
where the limits $ \dot{{\cal M}}(u,\theta) $'s satisfy that there exists a $ \Theta $-valued random vector $ \theta_* $ such that 
\begin{equation}\label{alternative 2}
\inf_{\theta: ||\theta-\theta_*||> \varepsilon}
||\dot{{\cal M}}(1,\theta)||
>0=||\dot{{\cal M}}(1,\theta_*)||, 
\quad \mbox{almost surely}, \quad \forall \varepsilon >0, 
\end{equation}
and that 
\begin{equation}\label{alternative 3}
\sup_{u \in (0,1)}||\dot{{\cal M}}(u,\theta_*)||>0, \quad \mbox{almost surely}. 
\end{equation}

\vskip 5pt
Assuming the conditions (\ref{null 1}), (\ref{null 2}), (\ref{alternative 1}) and (\ref{alternative 2}) is natural; 
see e.g.\ Theorems 5.7 and 5.9 of van der Vaart (1998). 
Let us explain how to check (\ref{alternative 3}) in the most typical form of alternatives in the change problems: 

\vskip 5pt
$ H_1' $: there exists a constant $ u_* \in (0,1) $ such that 
the true value is $ \theta_0 \in \Theta $ for $ u \in [0,u_*] $, 
and $ \theta_1 \in \Theta $ for $ u \in (u_*,1] $, 
where $ \theta_0 \not= \theta_1 $. 

\vskip 5pt
In many cases of ``ergodic models'', under $ H_1' $, 
the condition (\ref{alternative 1}) is satisfied with $ \dot{{\cal M}}(u,\theta) $ such that 
\[
\dot{{\cal M}}(u_*,\theta)=u_*\dot{M}_{\theta_0}(1,\theta) \quad \mbox{and} \quad 
\dot{{\cal M}}(1,\theta)=u_*\dot{M}_{\theta_0}(1,\theta)+(1-u_*)\dot{M}_{\theta_1}(1,\theta), 
\]
where not only $ \dot{M}_{\theta_0}(1,\theta) $'s but also $ \dot{M}_{\theta_1}(1,\theta) $'s are assumed to satisfy (\ref{null 2}) with trivial change of notation. 
To see that the condition (\ref{alternative 3}) is satisfied, notice that 
\[
\dot{{\cal M}}(u_*,\theta_*)=\dot{{\cal M}}(u_*,\theta_*)-u_*\dot{{\cal M}}(1,\theta_*)=u_*(1-u_*)
(\dot{M}_{\theta_0}(1,\theta_*)-\dot{M}_{\theta_1}(1,\theta_*)); 
\]
if this were zero with positive probability, then it should follow from $ \dot{{\cal M}}(1,\theta_*)=0 $ that 
$ \dot{M}_{\theta_0}(1,\theta_*) = \dot{M}_{\theta_1}(1,\theta_*) = 0 $ with positive probability, 
and this contradicts with (\ref{null 2}) and the assumption that $ \theta_0 \not=\theta_1 $. 
Therefore, we have 
\begin{eqnarray*}
\sup_{u \in (0,1)}||\dot{{\cal M}}(u,\theta_*)||&\geq& ||\dot{{\cal M}}(u_*,\theta_*)||
\\
&=&u_*(1-u_*)
||\dot{M}_{\theta_0}(\theta_*)-\dot{M}_{\theta_1}(\theta_*)|| > 0, 
\quad \mbox{almost surely}. 
\end{eqnarray*}
This positive value is closely related to the power of our test under $ H_1' $. 

Now, we prepare a lemma to prove the consistency of a sequence of $Z$-estimators. 
This lemma can be proved exactly in the same way as Theorems 5.7 and 5.9 of van der Vaart (1998), 
so the proof is omitted. 

\begin{lemma}\label{lemma consistency}
(i) Under [N], 
for any sequence of $ \Theta $-valued random vectors $ \widehat{\theta}_n $ 
such that $ ||\dot{\BM}_n(1,\widehat{\theta}_n)||=o_{P}(1) $, 
it holds that $ \widehat{\theta}_n \to^p \theta_0 $. 

(ii) Under [A], 
for any sequence of $ \Theta $-valued random vectors $ \widehat{\theta}_n $ 
such that $ ||\dot{\BM}_n(1,\widehat{\theta}_n)||=o_{P}(1) $, 
it holds that $ \widehat{\theta}_n \to^p \theta_* $. 
\end{lemma}

We are ready to state our main result of this section. 

\begin{theorem}\label{main 2}
Consider the above situation. 
Let $ \widehat{\theta}_n $ be any sequence of $ \Theta $-valued random 
vectors such that 
$ ||R_n \dot{\BM}_{n}(1,\widehat{\theta}_n)||=o_{P}(1) $ 
under $ H_0 $ and $ H_1 $. 
Let $ u \leadsto \widehat{V}_n(u) $ be any sequence of 
matrix valued random processes, which are regular except for $ u=0 $ almost 
surely, and it should be a uniformly consistent sequence of estimators 
for the non-negative definite matrix valued random process 
$ u \leadsto V(u,\theta_0) $ 
appearing below under $ H_0 $. 
Introduce the test statistic 
\[
{\cal T}_n=\sup_{u \in (0,1]}(R_n\dot{\BM}_n(u,\widehat{\theta}_n))^\top 
(u\widehat{V}_n(u)^{-1})R_n\dot{\BM}_n(u,\widehat{\theta}_n). 
\]

(i) Under [N], suppose that there exists a sequence of 
matrix valued random processes $ u \leadsto V_n(u,\theta_0) $ 
such that $ V_n(1,\theta_0) $'s are regular almost surely and that 
for any sequence of $ \Theta $-valued random vectors 
$ \widetilde{\theta}_n(u) $ indexed by $ u \in [0,1] $ 
satisfying $ \sup_{u \in [0,1]}||\widetilde{\theta}_n(u)-\theta_0 ||
\to^p 0 $, 
\begin{equation}\label{uniform convergence}
\sup_{u \in [0,1]}
||\ddot{\BM}_n(u,\widetilde{\theta}_n(u))+V_n(u,\theta_0)||\to^p 0. 
\end{equation}
Suppose also that 
\begin{equation}\label{sequential weak convergence}
(R_n \dot{\BM}_n(u,\theta_0),V_n(u,\theta_0))
\to^d ((u^{-1}V(u,\theta_0))^{1/2}B(u),V(u,\theta_0)),  
\quad \mbox{in } D[0,1], 
\end{equation}
where $ u \leadsto V(u,\theta_0) $ is 
a non-negative definite matrix valued random process 
such that $ V(1,\theta_0) $ is positive definite almost surely, 
and $ u \leadsto B(u) $ is 
a vector of independent standard Brownian motions; 
the value of the first vector of the limit 
in (\ref{sequential weak convergence}) 
at $ u=0 $ should be read as zero. 
(In general we do not assume that $ u \leadsto V(u,\theta_0) $ 
and $ u \leadsto B(u) $ are 
independent.) 

If $ \sup_{u \in [0,1]}||\widehat{V}_n(u)-V(u,\theta_0) ||\to^p 0 $, 
then it holds that 
\begin{equation}\label{main 2 conclusion}
{\cal T}_n\to^d  \sup_{u \in [0,1]}
||B(u)-u^{1/2}V(u,\theta_0)^{1/2}V(1,\theta_0)^{-1/2}B(1)||^2. 
\end{equation}
Therefore the test is asymptotically distribution free 
if $ V(u,\theta_0)=uV(1,\theta_0) $ for every $ u \in [0,1] $, 
because the limit in this case is reduced to $ \sup_{u \in [0,1]}||B^\circ(u)||^2 $ 
where $ u \leadsto B^\circ(u)=B(u)-uB(1) $ is a vector 
of independent standard Brwonian bridges. 
In the general case, if $ u \leadsto V(u,\theta_0) $ 
and $ u \leadsto B(u) $ are independent, then 
the limit in (\ref{main 2 conclusion}) is approximated by
\[
\sup_{u \in [0,1]}
||B(u)-u^{1/2}\widehat{V}_n(u)^{1/2}\widehat{V}_n(1)^{-1/2}B(1)||^2, 
\]
whose approximate distribution can be computed by some 
computer simulations for the standard Brownian motions $ u \leadsto B(u) $. 

(ii) Under [A], it holds for any random point $ \check{u} $ in $ (0,1) $ that 
\[
{\cal T}_n \geq 
\lambda(\check{u} R_n^2\widehat{V}_n(\check{u})^{-1})
\left\{ 
||\dot{{\cal M}}(\check{u},\theta_*)||^2 + o_{P}(1) \right\}, 
\]
where $ \lambda(A) $ denotes the smallest eigenvalue of the random matrix $ A $.  
Hence, if there exists a random point $ \check{u} $ in $ (0,1) $ such that 
$ ||\dot{{\cal M}}(\check{u},\theta_*)||>0 $ almost surely and that 
$ \lambda(R_n^{2}\widehat{V}_n(\check{u})^{-1}) $ tends to $ \infty $ 
in probability, then the test is consistent. 
\end{theorem}

{\bf Remark.} In the typical cases of ergodic models, the matrix $ V(u,\theta_0) $ is actually $ u I(\theta_0) $ where $ I(\theta_0) $ 
is the Fisher information matrix. Hence $ V(u,\theta_0)=uV(1,\theta_0) $ holds, and the reult is reduced to the standard case. 

\vskip 10pt
{\it Proof.} 
First let us prove (i). 
By Lemma \ref{lemma consistency} (i) 
we know that $ \widehat{\theta}_n $ is a consistent estimator for $ \theta_0 $ under $ H_0 $. 
So it follows from Theorem \ref{main theorem 1-1} that 
\begin{eqnarray*}
\lefteqn{
R_n\dot{\BM}_n(u,\widehat{\theta}_n)
}\\
&=&
R_n\dot{\BM}_n(u,\theta_0)
+\ddot{\BM}_n(u,\widetilde{\theta}_n(u))R_n(\widehat{\theta}_n-\theta_0)
\\
&=&
R_n\dot{\BM}_n(u,\theta_0)
-V_n(u,\theta_0)V_n(1,\theta_0)^{-1}R_n \dot{\BM}_n(1,\theta_0)
+e_n(u)
\\
&\to^d&
(u^{-1}V(u,\theta_0))^{1/2}(B(u)-u^{1/2}V(u,\theta_0)^{1/2}V(1,\theta_0)^{-1/2}
B(1)), \quad \mbox{in } D[0,1], 
\end{eqnarray*} 
where $ \widetilde{\theta}_n(u) $ is a random vector 
on the segment connecting $ \theta_0 $ and $ \widehat{\theta}_n $, 
and the reminder terms $ e_n(u) $ appearing above satisfy 
that $ \sup_{u \in [0,1]}||e_n(u)||\to^p 0 $. 
As a result the claim (i) follows 
from the continuous mapping theorem. 

The inequality in (ii) is proved as follows: 
\begin{eqnarray*}
{\cal T}_n &\geq& (R_n\dot{\BM}_n(\check{u},\widehat{\theta}_n))^\top
(\check{u}\widehat{V}_n(\check{u})^{-1}) R_n\dot{\BM}_n(\check{u},\widehat{\theta}_n) 
\\
&=& 
\dot{\BM}_n(\check{u},\widehat{\theta}_n)^\top
(\check{u} R_n^2\widehat{V}_n(\check{u})^{-1})\dot{\BM}_n(\check{u},\widehat{\theta}_n) 
\\
&\geq&
\lambda(\check{u}R_n^2\widehat{V}_n(\check{u})^{-1})||\dot{\BM}_n(\check{u},\widehat{\theta}_n)||^2
\\
&=&
\lambda(\check{u}R_n^2\widehat{V}_n(\check{u})^{-1})
\left\{ ||\dot{{\cal M}}(\check{u},\theta_*)||^2
+o_{P}(1) \right\}.  
\end{eqnarray*}
The proof is finished. 
\qed

\section{Example A: Ergodic diffusion process}\label{A}

Recall the description of Example A in Section \ref{introduction of examples}, 
where the first $ d_A $-components $ \alpha $ of the parameter $ \theta=(\alpha^\top,\beta^\top)^\top $ is involved in the drift coefficient, 
and the latter $ d_B $-components $ \beta $ is in the diffusion coefficient. 
Recalling also the definition of the rate matrix $ R_n $ there, let us 
consider the $ (d_A+d_B) $-dimensional random vectors $ \dot{M}_n(u,\theta) $ 
and the $ (d_A+d_B) \times (d_A+d_B) $-random matrices $ \ddot{M}_n(u,\theta) $ 
given as follows: 
\begin{eqnarray*}
\dot{M}_n(u,\theta)&=&(\dot{M}_n^A(u,\theta)^\top, 
\dot{M}_n^B(u,\theta)^\top)^\top, \\
\ddot{M}_n(u,\theta)&=&
\left(
\begin{array}{cc}
\ddot{M}_n^A(u,\theta) & \ddot{M}_n^C(u,\theta) \\
\ddot{M}_n^C(u,\theta)^\top  & \ddot{M}_n^B(u,\theta) 
\end{array}
\right). 
\end{eqnarray*}

Below, we will use the following notation: for a given constant $ p \geq 1 $ 
and a given sequence of positive constants $ r_n $, 
\begin{equation}\label{notation of high order moment}
\xi_n=o_{M(p)}(r_n^{-1}) \quad \Longleftrightarrow \quad r_nE[||\xi_n||^p] \to 0. 
\end{equation}
Notice that $ \xi_n=o_{M(1)}(r_n^{-1}) $ implies that 
$ \xi_{n}=o_{P}(r_n^{-1}) $. 

Under some regularlity 
conditions which are usually assumed in the asymptotic theory for ergodic diffusion process models, it is standard to show 
the followoing facts 
(see e.g. the appendix of Kessler (1997) for some techniques needed to prove them; 
see Nishiyama (2011), in Japanese, for the detailed proofs 
of the techniques that are omitted in Kessler's (1997) appendix):   
\[
\sup_{u \in [0,1]}
\left|\left|\dot{\BM}_n^A(u,\theta_0)
-\frac{1}{t_n^n}\sum_{k: t_{k-1}^n \leq u t_n^n}
\frac{\dot{S}(X_{t_{k-1}^n};\alpha_0)}{\sigma(X_{t_{k-1}^n};\beta_0)}
(W_{t_{k}^n}-W_{t_{k-1}^n})
\right|\right|=
o_{M(p)}((t_n^{n})^{-1/2}), 
\]
\[
\sup_{u \in [0,1]}
\left|\left|
\dot{\BM}_n^B(u,\theta_0)
-\frac{1}{n}\sum_{k: t_{k-1}^n \leq u t_n^n}
\frac{\dot{\sigma}(X_{t_{k-1}^n};\beta_0)}{\sigma(X_{t_{k-1}^n};\beta_0)}
\left\{\frac{|W_{t_{k}^n}-W_{t_{k-1}^n}|^2}{|t_k^n-t_{k-1}^n|} -1 \right\}
\right|\right|
=o_{M(p)}(n^{-1/2}), 
\]
\[
\sup_{u \in [0,1]}
\sup_{\theta \in \Theta}
\left|\left|\ddot{\BM}_n^A(u,\theta)-\frac{1}{t_n^n}
\sum_{k: t_{k-1}^n \leq u t_n^n}
H^A(X_{t_{k-1}^n};\theta_0,\theta) |t_{k}^n-t_{k-1}^n| 
\right|\right|
=o_{M(p)}((t_n^n)^{-1/2}), 
\]
\[
\sup_{u \in [0,1]}
\sup_{\theta \in \Theta}
\left|\left|
\ddot{\BM}_n^B(u,\theta)- \frac{1}{n}\sum_{k: t_{k-1}^n \leq u t_n^n}
H^B(X_{t_{k-1}^n};\theta_0,\theta) 
\right|\right|
=o_{M(p)}(n^{-1/2}), 
\]
\[
\sup_{u \in [0,1]}
\sup_{\theta \in \Theta}
||\ddot{\BM}_n^C(u,\theta)||
=o_{M(p)}(n^{-1/4}), 
\]
where 
\begin{eqnarray*}
H^A(x;\theta_0,\theta)&=&
\frac{\ddot{S}(x;\alpha)(S(x;\alpha_0)-S(x;\alpha))-\dot{S}(x;\alpha)\dot{S}(x;\alpha)^\top}
{\sigma(x;\beta)^2}, 
\\
H^B(x;\theta_0,\theta)&=&
\left\{\frac{\ddot{\sigma}(x;\beta)}{\sigma(x;\beta)^3} -3\frac{\dot{\sigma}(x;\beta)\dot{\sigma}(x;\beta)^\top}{\sigma(x;\beta)^4}
\right\}(\sigma(x;\beta_0)^2-\sigma(x;\beta)^2)
\\
& & 
-2 \frac{\dot{\sigma}(x;\beta)\dot{\sigma}(x;\beta)^\top}{\sigma(x;\beta)^2}. 
\end{eqnarray*}
The regularity conditions for the above claims depend on the constant 
$ p \geq 1 $ appearing in ``$ o_{M(p)}(r_n^{-1}) $'' which we need to have. 

\subsection{Moment convergence}

The assumption (\ref{main theorem 1-2 a condition}) can be checked 
by applying Burkholder-Davis-Gundy's inequality to the main part of 
$ R_n\dot{\BM}_n(\theta_0)=R_n\dot{\BM}_n(1,\theta_0) $. 
On the other hand, noting also $ t_n^n \leq n $, we can apply Remark 1 (ii) of Uchida and Yoshida (2012) to show that 
the assumption (\ref{main theorem 1-2 key condition}) 
for $ \ddot{\BM}_n(\theta)=\ddot{\BM}_n(1,\theta) $ is satisfied for 
\[
\ddot{M}(\theta)=
\left(
\begin{array}{cc}
\ddot{M}^A(\theta) & 0 \\
0 & \ddot{M}^B(\theta)
\end{array}
\right), 
\]
with 
\[
\ddot{M}^A(\theta)= \int_{I}H^A(x;\theta_0,\theta)\mu_{\theta_0}(dx) 
\quad \mbox{and} \quad \ddot{M}^B(\theta)= \int_{I}H^B(x;\theta_0,\theta)\mu_{\theta_0}(dx), 
\]
where $ \mu_{\theta_0} $ denotes the invariant distribution of $ X $ when the true value is $ \theta_0 $. 
In order to make the assumption [M1] or [M2] fulfilled, 
we have to introduce the parametric model for the drift and diffusion coefficients nicely. 
An example for which the assumption [M1] can be easily checked is $ S(\cdot;\alpha)=\alpha^\top a(\cdot) $ and $ \sigma(\cdot;\beta)=e^{\beta^\top b(\cdot)} $, 
where $ a(\cdot) $ and $ b(\cdot) $ are some vectors of known functions, assuming that $ b(\cdot) $ is bounded. 
The assumption [M2] would be satisfied in more general parametric models, because $ \ddot{M}(\theta) $'s are non-random in this example. 

\subsection{Change point problem}\label{Song Lee}

Under some standard conditions on the parametric family for the drift and diffusion coefficients in the context of 
ergodic diffusion processes, we can show that the condition (\ref{null 1}) under $ H_0 $ 
is satisfied with $ \dot{M}_{\theta_0}(1,\theta)= (\dot{M}_{\theta_0}^A(1,\theta)^\top, \dot{M}_{\theta_0}^B(1,\theta)^\top)^\top $, 
where 
\[
\dot{M}_{\theta_0}^A(1,\theta)=\int_{I}\frac{\dot{S}(x;\alpha)}{\sigma(x;\beta)}(S(x;\alpha_0)-S(x;\alpha))\mu_{\theta_0}(dx)
\]
and 
\[
\dot{M}_{\theta_0}^B(1,\theta)=\int_{I}\frac{\dot{\sigma}(x;\beta)}
{\sigma(x;\beta)^3}(\sigma(x;\beta_0)^2-\sigma(x;\beta)^2)\mu_{\theta_0}(dx), 
\]
and that the condition (\ref{alternative 1}) under $ H_1' $ is satisfied with 
\[
\dot{{\cal M}}(u,\theta)=(u \wedge u_*)\dot{M}_{\theta_0}(1,\theta)+ ((u-u_*) \vee 0)\dot{M}_{\theta_1}(1,\theta). 
\]
As stated there, the condition (\ref{alternative 3}) is automatically satisfied 
as soon as the natural conditions (\ref{null 2}) and (\ref{alternative 2}) are satisfied. 

Using the facts which we presented at the beginning of this section and 
the usual martingale central limit theorem, we can see that 
the condition (\ref{uniform convergence}) and (\ref{sequential weak convergence}) 
hold for 
\[
V_n(u,\theta_0)=
\left(
\begin{array}{cc}
V_n^A(u,\theta_0) & 0 \\
0 & V_n^B(u,\theta_0)
\end{array}
\right), 
\]
where 
\begin{eqnarray*}
V_n^{A}(u,\theta_0)&=&\frac{1}{t_n^n}
\sum_{k: t_{k-1}^n \leq u t_n^n}
\frac{\dot{S}(X_{t_{k-1}^n};\alpha_0)
\dot{S}(X_{t_{k-1}^n};\alpha_0)^\top}
{\sigma(X_{t_{k-1}^n};\beta_0)^2}
|t_{k}^n-t_{k-1}^n|,\\
V_n^{B}(u,\theta_0)&=&\frac{2}{n}
\sum_{k: t_{k-1}^n \leq u t_n^n}
\frac{\dot{\sigma}(X_{t_{k-1}^n};\beta_0)
\dot{\sigma}(X_{t_{k-1}^n};\beta_0)^\top}
{\sigma(X_{t_{k-1}^n};\beta_0)^2}. 
\end{eqnarray*}
The limit of $ V_n(u,\theta_0) $ is 
$ V(u,\theta_0)=u I_{\theta_0}(\theta_0) $, 
where 
\[
I_{\theta_0}(\theta)=
\left(
\begin{array}{cc}
I_{\theta_0}^A(\theta) & 0 \\
0 & I_{\theta_0}^B(\theta)
\end{array}
\right)
\]
with 
\[
I_{\theta_0}^A(\theta)=\int_I\frac{\dot{S}(x;\alpha)\dot{S}(x;\alpha)^\top}{\sigma(x;\beta)^2}\mu_{\theta_0}(dx),
\quad 
I_{\theta_0}^B(\theta)=2\int_I\frac{\dot{\sigma}(x;\beta)\dot{\sigma}(x;\beta)^\top}{\sigma(x;\beta)^2}\mu_{\theta_0}(dx). 
\]
We suppose that $ I_{\theta_0}(\theta) $'s are positive definite. 

As a consistent estimator $ \widehat{V}_n(u) $ for $ V(u,\theta_0) $, 
we introduce 
\[
\widehat{V}_n(u)=
\left(
\begin{array}{cc}
u \widehat{I}_n^A & 0 \\
0 & u \widehat{I}_n^B 
\end{array}
\right),
\]
where 
\[
\widehat{I}_n^A= \frac{1}{n}\sum_{k=1}^n \frac{\dot{S}(X_{t_{k-1}^n};\widehat{\alpha}_n)\dot{S}(X_{t_{k-1}^n};\widehat{\alpha}_n)^\top}
{\sigma(X_{t_{k-1}^n};\widehat{\beta}_n)^2}, 
\quad 
\widehat{I}_n^B=
\frac{2}{n}\sum_{k=1}^n \frac{\dot{\sigma}(X_{t_{k-1}^n};\widehat{\beta}_n)\dot{\sigma}(X_{t_{k-1}^n};\widehat{\beta}_n)^\top}
{\sigma(X_{t_{k-1}^n};\widehat{\beta}_n)^2}. 
\]
Since $ V(u,\theta_0)=uV(1,\theta_0) $ in this example, 
the limit of 
our test statistic $ {\cal T}_n $ is $ \sup_{u \in [0,1]}||B^\circ(u)||^2 $ where $ u \leadsto B^\circ(u) $ is a vector of standard Brownian bridges, 
so the test is asymptotically distribution free. 

Finally, it is clear that under $ H_1' $, $ \lambda(R_n^2\widehat{V}_n(u_*)^{-1}) $ tends to $ \infty $ in probability 
since the matrix 
$ u_*I_{\theta_0}(\theta_*) + (1-u_*)I_{\theta_1}(\theta_*) $, which is the limit of $ \widehat{V}_n(u_*) $, is positive definite. 
Thus the test is consistent. 

\subsection{Numerical study for change point problem}

In this section, as well as Section \ref{numerical B} for Example B, we observe finite sample performance of our test statistic 
through numerical experiments. Here, we adopt the Ornstein-Uhlenbeck process starting from $x_0=0$ 
for the true (data-generating) process: 
\begin{equation}
X_{t}=x_0-\int_{0}^{t}\alpha X_{s}ds+\beta W_{t}, \quad t\in [0,T].
\label{DGP}
\end{equation}
For simplicity, we shall treat the equidistant sampling case, that is, $ \Delta_{n}=|t^{n}_{k}-t^{n}_{k-1}|$ 
for every $ k=1,...,n $. 

We are going to observe the trajectory of the process (\ref{DGP}) 
for different time horizons $t_n^n=T$, and the number $n$  of observations for each trajectory is such that 
$t_n^n=n^{1/3}$, so $\Delta_n= n^{-2/3}$. 
For this process (\ref{DGP}) the estimators for the parameters $\alpha$ and $\beta$ and 
the estimator of the information matrix can be explicitly calculated, and thus the test statistic can be easily computed. 
For any fixed level $\varepsilon>0$ the critical value $c_\varepsilon$ is given by 
\[
P\left(\sup_{u \in [0,1]}\sum_{i=1}^{d_A +d_B}|B^{\circ,(i)}(u)|^2>c_\varepsilon\right)=\varepsilon. 
\]
Table 1 of Lee {\em et al.} (2003) gives a table of the critical values for the significance levels 
$\varepsilon=0.01, 0.05, 0.10$ and for different values of the dimension $d=d_A+d_B$ computed by 
Monte Carlo simulation for the limit distribution. 
Throughout we take the significance level to be $\varepsilon=0.05$. For two parameters ($d=2$) the critical value is $c_\varepsilon=2.408$. 
Regarding the null hypothesis we generate $M=10^4$ trajectory of (\ref{DGP})  and we evaluate the empirical size. 
The results are reported in Table \ref{table1}. We observe that: 
the empirical size gains along with increasing terminal time $T=t^{n}_{n}$, attaining at 0.05, but also for small terminal $T$. 
In the second example reported in Table \ref{table1}, the values of the parameter are the maximum likelihood estimate for the mostly federal funds data 1963-1998 
in A\"it-Sahalia (1999). 

\begin{table}[h]
\begin{center}
\begin{tabular}{r|rrrrrrr}
\hline
 $T$& $5$& $10$ &15 &$20$ & $25$\\
$n$ & {\tiny 125}& {\tiny 1000} & {\tiny 3375} & {\tiny 8000}& {\tiny 15625}\\
\hline
&&&&&\\
$\alpha_0=1$, $\beta_0=1$& 0.044 &  0.054 & 0.050  &  0.052  & 0.053\\
&&&&&\\
\hline
&&&&&\\
$\alpha_0=0.25$, $\beta_0=0.02$&  0.047& 0.061 & 0.058  & 0.064 & 0.054\\
&&&&&\\
\hline 
\end{tabular}
\caption{Empirical size based on $M=10^4$ independent statistics, for different time horizons.} 
\label{table1}
\end{center}
\end{table}

Regarding the alternative hypothesis we study the behavior of the test statistic 
in three different situations and for different change point $ u_*T $ of the parameters, as follows: 
\begin{itemize}
\item The drift coefficient changes from $\alpha_0$ to $\alpha_1$, but the diffusion coefficient does not change. 
\item The drift coefficient does not change, but the diffusion coefficient changes from $\beta_0$ to $\beta_1$. 
\item Both coefficients change. 
\end{itemize}  
For each of the above scenarios we consider the following change points, $u_*=\frac12, \frac34, \frac{9}{10}$. 

The first scenario is the worst case for the diffusion, and in order to detect a change in the drift 
we have to observe the process as long as possible. 
Table \ref{table2} shows empirical power for different terminal times $T$ and different change points $u_*T $.  
The values of the parameters are 
$\alpha_0=0.25$, $\alpha_1=0.50$ and $\beta=0.02$. The last does not vary. 
We simulate  $10^4 $ independent copies of a trajectory of (\ref{DGP}) to obtain different values of 
$ {\cal T}_n$. 
The power increase as $T$ increase and the performance is better when we can observe the process after the change for long time 
(the case $u_*=\frac12$). 
In such a case the power of the test is reasonable. 
In the worst case $ u_*=\frac{9}{10} $, the test is not able to detect the change in the drift coefficient. 
\begin{table}[h]
\begin{center}
\begin{tabular}{r|rrrrrrr}
\hline
$T$& $5$& $10$ &15 &$20$ & $25$\\
$n$ & {\tiny125}& {\tiny1000} & {\tiny 3375} & {\tiny 8000}& {\tiny 15625}\\
\hline
&&&&&\\
$u_*=\frac12$ & 0.31 & 0.52  & 0.73  &  0.79  & 0.88\\
&&&&&\\
\hline
&&&&&\\
$u_*=\frac34$&  0.12& 0.17 & 0.23  & 0.26 & 0.35\\
&&&&&\\
\hline 
&&&&&\\
$u_*=\frac{9}{10}$&  0.05& 0.07 & 0.08  & 0.08 & 0.09\\
&&&&&\\
\hline 
\end{tabular}
\caption{Empirical power based on $M=10^4$ independent statistics. 
Here the significance level is $0.05$. 
The values of the parameter are $\alpha_0=0.25$, $\alpha_1=0.50$ and $\beta=0.02$ (it does not vary).}
\label{table2}
\end{center}
\end{table}

Table \ref{table3} reports the results for simulation when only the drift changes, but the change is bigger. 
With $\alpha_0=0.25$, $\alpha_1=1.25$ and $\beta=0.02$, the power increases not only for $u_*=\frac12 $ but also for $u_*=\frac34 $. 
This was expected, but as in the previous example the performance of the test is not good when the chance of the parameter occur 
at the end of the observation window. 
\begin{table}[h]
\begin{center}
\begin{tabular}{r|rrrrrrr}
\hline
$T$& $5$& $10$ &15 &$20$ & $25$\\
$n$ & {\tiny125}& {\tiny1000} & {\tiny 3375} & {\tiny 8000}& {\tiny 15625}\\
\hline
&&&&&\\
$u_*=\frac12$ & 0.35 & 0.60  & 0.78 &  0.88 & 0.94\\
&&&&&\\
\hline
&&&&&\\
$u_*=\frac34$&  0.13 & 0.20 & 0.28  & 0.31 & 0.38\\
&&&&&\\
\hline 
&&&&&\\
$u_*=\frac{9}{10}$&  0.06 & 0.08 & 0.09 & 0.11 & 0.11\\
&&&&&\\
\hline 
\end{tabular}
\caption{Empirical power based on $M=10^4$ independent statistics. 
Here the significance level is $0.05$. The values of the parameter are $\alpha_0=0.25$, $\alpha_1=1.25$ 
and $\beta=0.02$ (it does not vary).}
\label{table3}
\end{center}
\end{table}

Table \ref{table4} shows the empirical power for different terminal times $T$ and different change points $u_*$ in the second scenario: 
$\alpha=0.25$ does not vary, but 
$\beta_0=0.02$ changes and becomes $\beta_1=0.03$.  
The power of the test is very good, but this is not surprising because a change in the diffusion coefficient can be easily detected. 
The situation reported in Table \ref{table5} is the same: also for very small change in the diffusion coefficient, 
the performance of the test is very good with the empirical power that reaches the value 1 also for small $T$. 
We do not report the results for the third scenario where the drift changes at the same instant of the diffusion, 
because the performance of the test is the same as in the second scenario. 
This is not surprising and is due to the different rates of convergence of the estimators of the two parameters.  
\begin{table}[h]
\begin{center}
\begin{tabular}{r|rrrrrrr}
\hline
$T$& $5$& $10$ &15 &$20$ & $25$\\
$n$ & {\tiny125}& {\tiny1000} & {\tiny 3375} & {\tiny 8000}& {\tiny 15625}\\
\hline
&&&&&\\
$u_*=\frac12$ & 0.99& 1  & 1  &  1  & 1 \\
&&&&&\\
\hline
&&&&&\\
$u_*=\frac34$&  0.86& 1 & 1  & 1 & 1 \\
&&&&&\\
\hline 
&&&&&\\
$u_*=\frac{9}{10}$&  0.36& 0.99 & 1  & 1 & 1\\
&&&&&\\
\hline 
\end{tabular}
\caption{Empirical power based on $M=10^4$ independent statistics. 
Here the significance level is $0.05$. The values of the parameter are $\alpha=0.25$ (it does not vary), 
$\beta_0=0.02$ and $\beta_1=0.03$. }
\label{table4}
\end{center}
\end{table}

\begin{table}[h]
\begin{center}
\begin{tabular}{r|rrrrrrr}
\hline
$T$& $5$& $10$ &15 &$20$ & $25$\\
$n$& {\tiny125}& {\tiny1000} & {\tiny 3375} & {\tiny 8000}& {\tiny 15625}\\
\hline
&&&&&\\
$u_*=\frac12$ & 0.87& 1  & 1  &  1  & 1\\
&&&&&\\
\hline
&&&&&\\
$u_*=\frac34$&  0.52& 0.99 & 1  & 1 & 1\\
&&&&&\\
\hline 
&&&&&\\
$u_*=\frac{9}{10}$&  0.14& 0.62 & 0.99  & 1 & 1\\
&&&&&\\
\hline 
\end{tabular}
\caption{Empirical power based on $M=10^4$ independent statistics. Here the significance level is $0.05$. 
The values of the parameter are $\alpha=0.25$ (it does not vary), 
$\beta_0=0.020$ and $\beta_1=0.025$. }
\label{table5}
\end{center}
\end{table}

\section{Example B: Volatility of diffusion process}\label{B}

Recall the description of Example B in Section \ref{introduction of examples}. 
An interesting point of this example is that the limit of 
$ -\ddot{\BM}_n(u,\widetilde{\theta}_n(u)) $ is random and 
depend on $ u \in [0,1] $ in a complex way. 

Let a constant $ p \geq 1 $ be given, and 
recall the notation (\ref{notation of high order moment}). 
Under some regularity conditions, it holds that 
\[
\sup_{u \in [0,1]}
\left|\left|
\dot{\BM}_n(u,\theta_0)
-\frac{1}{n}\sum_{k: t_{k-1}^n \leq u t_n^n}
\frac{\dot{\sigma}(X_{t_{k-1}^n};\theta_0)}{\sigma(X_{t_{k-1}^n};\theta_0)}
\left\{\frac{|W_{t_{k}^n}-W_{t_{k-1}^n}|^2}{|t_k^n-t_{k-1}^n|} -1 \right\}
\right|\right|
=o_{M(p)}(n^{-1/2}),
\] 
\[
\sup_{u \in [0,1]}
\sup_{\theta \in \Theta}
\left|\left|
\ddot{\BM}_n(u,\theta)- \frac{1}{n}\sum_{k: t_{k-1}^n \leq u t_n^n}
H(X_{t_{k-1}^n};\theta_0,\theta) 
\right|\right|
=o_{M(p)}(n^{-1/2}), 
\]
where 
\begin{eqnarray*}
H(x;\theta_0,\theta)&=&
\left\{\frac{\ddot{\sigma}(x;\theta)}{\sigma(x;\theta)^3} -3\frac{\dot{\sigma}(x;\theta)\dot{\sigma}(x;\theta)^\top}{\sigma(x;\theta)^4}
\right\}(\sigma(x;\theta_0)^2-\sigma(x;\theta)^2)
\\
& & 
-2 \frac{\dot{\sigma}(x;\theta)\dot{\sigma}(x;\theta)^\top}{\sigma(x;\theta)^2}. 
\end{eqnarray*}
The regularity conditions for the above claims 
depend on the constant $ p \geq 1 $ which we need to have. 
Moreover, 
under some standard conditions, 
it holds that for any sequence of random vectors $ \widetilde{\theta}_n(u) $ indexed by $ u \in [0,1] $ 
such that $ \sup_{u \in [0,1]}||\widetilde{\theta}_n(u)-\theta_0|| \to^p 0 $, 
\[
\sup_{u \in [0,1]}||\ddot{\BM}_n(u,\widetilde{\theta}_n(u))
+V_n(u,\theta_0)|| \to^p 0, 
\]
where
\[
V_n(u,\theta_0)=\frac{2}{n}\sum_{k: t_{k-1}^n \leq u t_n^n}
\frac{\dot{\sigma}(X_{t_{k-1}^n};\theta_0)\dot{\sigma}(X_{t_{k-1}^n};\theta_0)^\top}{\sigma(X_{t_{k-1}^n};\theta_0)^2}, \quad \forall u\in [0,1]. 
\]
Also, it follows from the well known theory of martingales that 
\[
(\sqrt{n}\dot{\BM}_n(u,\theta_0),V_n(u,\theta_0))
\to^d ((u^{-1}V(u,\theta_0))^{1/2}B(u),V(u,\theta_0))
\quad \mbox{in }D[0,1], 
\]
where $ u \leadsto B(u) $ is a vector of independent standard Brownian motions which is independent of the 
matrix valued random process $ u \leadsto V(u,\theta_0) $ given by 
\[
V(u,\theta_0)=2\int_{0}^{uT}
\frac{\dot{\sigma}(X_s;\theta_0)\dot{\sigma}(X_s;\theta_0)^\top}{\sigma(X_s;\theta_0)^2}ds, \quad \forall u \in [0,1]. 
\]

\subsection{Moment convergence}

Due to the above facts (for $ u=1 $), Theorem \ref{main theorem 1-1} yields that 
for any consistent estimator $ \widehat{\theta}_n $ for $ \theta_0 $ satisfying 
$ ||\dot{\BM}_n(\widehat{\theta}_n)||=o_{P}(n^{-1/2}) $ we have 
$ \sqrt{n}(\widehat{\theta}_n-\theta_0)\to^d V(\theta_0)^{-1/2}Z $, 
where $ Z $ is a standard Gaussian random vector which is independent of $ V(\theta_0)=V(1,\theta_0) $. 

Next let us apply Theorem \ref{main theorem 1-2}. 
The assumption (\ref{main theorem 1-2 a condition}) 
for $ \sqrt{n}\dot{\BM}_n(\theta_0)=\sqrt{n}\dot{\BM}_n(1,\theta_0) $ 
can be checked by using Burkholder-Davis-Gundy's inequality. 
In the case of this example, 
checking that the assumption (\ref{main theorem 1-2 key condition}) 
for $ \ddot{\BM}_n(\theta)=\ddot{\BM}_n(1,\theta) $ 
is satisfied with 
\[
\ddot{M}(\theta)= \int_0^TH(X_t;\theta_0,\theta)dt
\]
is easy. 
In order to make the assumption [M1] or [M2] fulfilled, 
we again have to introduce the parametric model for the diffusion coefficients nicely. 
An example for which the former assumption in [M1] can be easily checked 
is $ \sigma(\cdot;\theta)=e^{\theta^\top g(\cdot)} $, 
where $ g(\cdot) $ are some vectors of known, bounded functions. 
The latter assumption in [M1] is then reduced to 
\[
E\left[\left|\left| \left(\int_0^T g(X_t)g(X_t)^\top dt \right)^{-1}\right|\right|^{pb/\gamma} \right] < \infty,  
\]
for which we can give a clear sufficient condition for the function $ g $ at least in the one-dimensional case 
(for example, just assume $ |g(\cdot)|^2 \geq c $ for a constant $ c>0 $). 

\subsection{Change point problem}\label{B2}

Under some standard conditions on the parametric family for the diffusion coefficient, 
we can show that the condition (\ref{null 1}) under $ H_0 $ is satisfied with 
\[
\dot{M}_{\theta_0}(1,\theta)=\int_{0}^{uT}\frac{\dot{\sigma}(X_t;\theta)}
{\sigma(X_t;\theta)^3}(\sigma(X_t;\theta_0)^2-\sigma(X_t;\theta)^2)dt. 
\]
Under $ H_1' $, we have 
\[
\sup_{u \in [0,1]}||\dot{\BM}_n(u,\theta)-\dot{{\cal M}}(u,\theta)|| \to^p 0, 
\]
where 
\begin{eqnarray*}
\dot{{\cal M}}(u,\theta)&=&\int_{0}^{(u \wedge u_*)T}\frac{\dot{\sigma}(X_t;\theta)}
{\sigma(X_t;\theta)^3}(\sigma(X_t;\theta_0)^2-\sigma(X_t;\theta)^2)dt
\\
& & 
+1\{ u > u_* \} \int_{u_*T}^{uT}\frac{\dot{\sigma}(X_t;\theta)}
{\sigma(X_t;\theta)^3}(\sigma(X_t;\theta_1)^2-\sigma(X_t;\theta)^2)dt. 
\end{eqnarray*}
We can give a set of sufficient conditions for (\ref{alternative 3}) as follows. 
Suppose that the Lebesgue measure of the random set 
$ T_0=\{t \in [0,T]: \dot{\sigma}(X_t;\theta_*)/\sigma(X_t;\theta_*)^3=0 \} $ 
is zero almost surely, which is true in many concrete models. 
In this case, replace the values 
$ \dot{\sigma}(X_t;\theta_*)/\sigma(X_t;\theta_*)^3 $ in the definition of 
$ \dot{{\cal M}}(u,\theta_*) $ on the set $ T_0 $ by $1$ to construct 
$ \dot{{\cal M}}(u,\theta_*)^\star $ which equals with 
the original $ \dot{{\cal M}}(u,\theta_*) $ for all $ u \in [0,1] $, 
almost surely. 
If we further assume that for any non-empty interval $ J \subset I $ 
\[
\sigma(x;\theta)=\sigma(x;\theta'), \quad \forall x \in J 
\quad \Longleftrightarrow \quad \theta=\theta'
\]
and that each of the random sets $ J_0=\{ X_t: t \in [0,u_*T] \} $ 
and $ J_1=\{ X_t: t \in (u_*T,T] \} $ includes a non-empty set almost surely, 
then it follows from the assumption $ \theta_0 \not= \theta_1 $ that 
$ \frac{d}{du}\dot{{\cal M}}(u,\theta_*)^\star \not\equiv 0 $ 
almost surely. 
Thus we have 
\[
\sup_{u \in (0,1)}||\dot{{\cal M}}(u,\theta_*)||
=\sup_{u \in (0,1)}||\dot{{\cal M}}(u,\theta_*)^\star||>0, \quad 
\mbox{almost surely}. 
\]

Now, consider the matrices 
\[
V(u,\theta)=2\int_0^{uT}\frac{\dot{\sigma}(X_t;\theta)\dot{\sigma}(X_t;\theta)^\top}{\sigma(X_t;\theta)^2}dt, \quad \forall u \in [0,1]. 
\]
Let us assume that for every $ \theta \in \Theta$ there exists a set $ J_\theta \subset I $ such that 
the Lebesgue measure of $ J_\theta^c $ is zero and that 
the matrices $ \dot{\sigma}(x;\theta)\dot{\sigma}(x;\theta)^\top/\sigma(x;\theta) $ are 
positive definite for $ x \in J_\theta $, which is a standard assumption. 
In this case, if the claim that the Lebesgue measure of the set $ \{ t \in [0,T]; X_t(\omega) \in J \} $ 
is positive for any set $ J \subset I $ such that the Lebesgue measure of $ J^c $ is zero holds for almost all $ \omega $, 
then $ V(u,\theta_0) $'s and $ V(u,\theta_*) $'s for $ u \in (0,1] $ are positive definite 
almost surely under $ H_0 $ and $ H_1' $, respectively.   

As we saw at the beginning of this section, 
the conditions (\ref{uniform convergence}) 
and (\ref{sequential weak convergence}) under $ H_0 $ 
are satisfied. 
As a consistent estimator $ \widehat{V}_n(u) $ for $ V(u,\theta_0) $ is given by 
\[
\widehat{V}_n(u)=
\frac{2}{n}\sum_{k: t_{k-1}^n \leq u t_n^n}\frac{\dot{\sigma}(X_{t_{k-1}^n};\widehat{\theta}_n)\dot{\sigma}(X_{t_{k-1}^n};\widehat{\theta}_n)^\top}
{\sigma(X_{t_{k-1}^n};\widehat{\theta}_n)^2}, \quad \forall u \in [0,1]. 
\]
Our test in this example is not asymptotically distribution free. 

Finally, it is clear that under $H_1'$, $ n \lambda(\widehat{V}_n(\check{u})^{-1}) $ tends to $ \infty $ in probability, 
because it follows from what we have assumed that 
$ \lambda(\widehat{V}_n(\check{u})^{-1}) \to^p 
\lambda(V(\check{u},\theta_*)^{-1}) $ and the limit is positive almost surely. 
Thus the test is consistent. 

\subsection{Numerical study for change point problem}\label{numerical B}

The data-generating process is the following: 
\[
X_t=4 - \int_0^{t}(X_s-4)ds +\int_{0}^{t}\exp\left(\theta \frac{X_s^2}{1+X_s^2}\right)dW_s, \quad t \in [0,1], 
\]
where the drift coefficient $ S(x)=-(x-4) $ is treated as a nuisance function. 
Suppose that we observe $ M=10^3 $ independent copies of this process at the equidistant time grid 
$ t_k^n=\frac{k}{n} $, $ k=0,1,...,n $. 
We compute the critical value of the test based on the approximation of the limit distribution 
\begin{equation}\label{5.3 limit}
\sup_{u \in[0,1]}|B(u)-u^{1/2}V(u,\theta_0)^{1/2}V(1,\theta_0)^{-1/2}B(1)|^2
\end{equation}
obtained by replacing 
\[
V(u,\theta_0) =2 \int_{0}^{u}\left|\frac{X_s^2}{1+X_s^2}\right|^2ds
\]
by the natural estimator 
\[
\widehat{V}_n(u)=\frac{2}{n}\sum_{k=1}^{[un]}\left|\frac{X_{t_{k-1}^n}^2}{1+X_{t_{k-1}^n}^2}\right|^2
\]
and doting $ 10^3 $ times Monte Carlo simulation for 
the standard Brownian motion $ u \leadsto B(u) $. 

The empirical size under $ H_0 $ is reported in Table \ref{table6}, where 
the true value of the parameter is set as $ \theta_0=1.0 $ or $ 1.5 $. 
We see that the convergence to the approximate distribution of (\ref{5.3 limit}) is not perfectly good, 
but it is reasonable even for the cases where $ n $ is small. 
\begin{table}[h]
\begin{center}
\begin{tabular}{r|rrrr}
\hline
$n$&$20$&$40$&$100$&$200$\\
\hline
&&&&\\
$\theta_0=1.0$&0.026&0.024&0.042& 0.040 \\
&&&&\\
\hline
&&&&\\
$\theta_0=1.5$&0.026&0.023&0.037&0.034 \\
&&&&\\
\hline 
\end{tabular}
\caption{Empirical size based on $M=10^3$ independent statistics. Here the significance level is $0.05$. 
The value of the parameter is $\theta_0=1.0 $ or $ 1.5 $.} 
\label{table6}
\end{center}
\end{table}

The empirical power under $ H_1' $ is reported in Table \ref{table7}, where 
the true values of the parameter change from $ \theta_0=1.0 $ to $ \theta_1=1.5 $ at time point $ u_* =\frac{1}{2}, \frac{3}{4} $ or $ \frac{9}{10}$. 
\begin{table}[h]
\begin{center}
\begin{tabular}{r|rrrr}
\hline
$n$&$20$&$40$&$100$&$200$\\
\hline
&&&&\\
$u_*=\frac{1}{2}$&0.067&0.331&0.755&0.946\\
&&&&\\
\hline
&&&&\\
$u_*=\frac{3}{4}$&0.125&0.255&0.630&0.873\\
&&&&\\
\hline
&&&&\\
$u_*=\frac{9}{10}$&0.048&0.117&0.275&0.462 \\
&&&&\\
\hline 
\end{tabular}
\caption{Empirical power based on $M=10^3$ independent statistics. Here the significance level is $0.05$. 
The values of the parameter change from $\theta_0=1.0 $ to $ \theta_1=1.5 $ at time $ u_*=\frac{1}{2}, \frac{3}{4} $ or $ \frac{9}{10}$.} 
\label{table7}
\end{center}
\end{table}

\section{Example C: Cox's regression model}\label{C}

Recall the description of Example C in Section \ref{introduction of examples}. 
Since all the arguments are similar to those in Section \ref{B}, 
we state only the key points in the discussion on the change point problem. 

Introducing the notations 
\begin{eqnarray*}
S_t^{n,0}(\theta)&=&\sum_{k=1}^n e^{\theta Z_t^k} Y_t^k, 
\\
S_t^{n,1}(\theta)&=&\sum_{k=1}^n Z_t^k e^{\theta Z_t^k} Y_t^k,
\\
S_t^{n,2}(\theta)&=&\sum_{k=1}^n (Z_t^{k})^\top Z_t^k e^{\theta Z_t^k} Y_t^k, 
\end{eqnarray*}
we suppose that 
\[
\sup_{\theta \in \Theta}\sup_{t \in [0,T]}\left|\left|\frac{1}{n}S_t^{n,l}(\theta)-{\cal S}_t^{l}(\theta)\right|\right| \to^p 0, 
\quad l=0,1,2, 
\]
where the limits $ t \leadsto {\cal S}_t^{l} $ are some stochastic processes 
(c.f. Andersen and Gill (1982) who assumed that $ {\cal S}^{l} $'s  are not random). 

Then, some arguments similar to Section \ref{B2} are possible for 
\begin{eqnarray*}
\dot{M}_{\theta_0}(1,\theta)&=&
\int_0^T \left(\frac{{\cal S}_t^{1}(\theta_0)}{{\cal S}_t^{0}(\theta_0)}-\frac{{\cal S}_t^{1}(\theta)}{{\cal S}_t^{0}(\theta)}\right)
{{\cal S}_t^{0}(\theta_0)}\alpha(t)dt, 
\\
\dot{{\cal M}}(u,\theta)&=&
\int_0^{(u \wedge u_*)T} \left(\frac{{\cal S}_t^{1}(\theta_0)}{{\cal S}_t^{0}(\theta_0)}-\frac{{\cal S}_t^{1}(\theta)}{{\cal S}_t^{0}(\theta)}\right)
{{\cal S}_t^{0}(\theta_0)}\alpha(t)dt
\\
& & 
+ 1\{ u > u_* \}
\int_{u_*T}^{uT}
\left(\frac{{\cal S}_t^{1}(\theta_1)}{{\cal S}_t^{0}(\theta_1)}-\frac{{\cal S}_t^{1}(\theta)}{{\cal S}_t^{0}(\theta)}\right)
{{\cal S}_t^{0}(\theta_1)}\alpha(t)dt, 
\\
V_n(u,\theta_0)&=&
\frac{1}{n}\int_{0}^{uT}\frac{S_t^{n,0}(\theta_0)S_t^{n,2}(\theta_0)-S_t^{n,1}(\theta_0)S_t^{n,1}(\theta_0)^\top}
{S_t^{n,0}(\theta_0)}\alpha(t)dt, 
\\
V(u,\theta)&=&\int_{0}^{(u \wedge u_*)T}
\frac{{\cal S}_t^{0}(\theta){\cal S}_t^{2}(\theta) - {\cal S}_t^{1}(\theta){\cal S}_t^{1}(\theta)^\top}
{{\cal S}_t^{0}(\theta)^2}{\cal S}_t^{0}(\theta_0)\alpha(t)dt
\\
& & 
+ 1\{ u > u_* \}
\int_{u_*T}^{uT}
\frac{{\cal S}_t^{0}(\theta){\cal S}_t^{2}(\theta) - {\cal S}_t^{1}(\theta){\cal S}_t^{1}(\theta)^\top}
{{\cal S}_t^{0}(\theta)^2}{\cal S}_t^{0}(\theta_1)\alpha(t)dt, 
\\
\widehat{V}_n(u)&=&
\frac{1}{n}
\sum_{k=1}^n
\int_{0}^{uT}
\frac{S_t^{n,0}(\widehat{\theta}_n)S_t^{n,2}(\widehat{\theta}_n) - S_t^{n,1}(\widehat{\theta}_n)S_t^{n,1}(\widehat{\theta}_n)^\top}
{S_t^{n,0}(\widehat{\theta}_n)^2}dN_t^k. 
\end{eqnarray*}

\vskip 20pt
\par\noindent
{\bf Acknowledgements.} 
This work was supported by Prin09 grant (I.N.) and 
by Grant-in-Aid for Scientific Research (C), 24540152, 
from Japan Society for the Promotion of Science (Y.N.).

\end{document}